\documentclass{article}

\def \R {{\mathbf {R}}}

\title{ Gaussian  factors, spectra, and $P$-entropy}
\author{ Valery V. Ryzhikov}
\date{}

\usepackage[T1]{fontenc} 
\usepackage[cp1251]{inputenc}  
\usepackage[russian]{babel}
\begin{document}

\begin{center}   
{ \Large \bf

Gaussian  factors, spectra, and $P$-entropy

\vspace{5mm}
Valery V. Ryzhikov}

\vspace{4mm} 
{\large There is  an infinite set (of continual cardinality) of pairwise disjoint Gaussian \\ automorphisms with  spectrally isomorphic  even factors with Lebesgue spectrum. }

\end{center}  


 

\Large
\section{Gaussian automorphisms with completely positive P-entropy}
Among  the transformations  $T$  preserving  infinite measure  there are a lot  of    $T$ with  singular spectrum   whose convolution powers are equivalent to Lebesgue \cite{24}.
The corresponding Gaussian automorphisms $G(T)$ and Poisson suspensions $P(T)$ (see \cite{KSF},\cite{NP}, \cite{PR}) over such $T$  have zero  entropy.   The metric properties  of $G(T)$  and  $P(T)$ can be very  different: for example, Gaussian automorphisms are infinitelly divisible, hence,  $G(T)$  have large centralizers, while  $P(T)$ have  the trivial ones (see \cite{06}, \cite{PR}). The even factors $G_{ev}(T)$ (\cite{NP}) of such Gaussian  automorphisms $G(T)$ are spectrally isomorphic:  their spectrum is  Lebesgue of  infinite multiplicity.   What non-spectral invariants of these factors  can distinguish them?  We use entropy invariants to show that 
\it there is an uncountable family of  orthogonal operators $T$ such that the corresponding pairwise 
non-isomorphic Gaussian even factors $G_{ev}(T)$ have the same spectrum but different entropy properties.  \rm

\vspace{3mm}
\bf $P$-entropy. \rm
Let $P=\{P_j\}$ be a sequence of finite subsets in a countable infinite group. For a Probability-preserving action $T$ of this group we define
$$h_j(T,\xi)=\frac 1 {|P_j|} H\left(\bigvee_{p\in P_j}T_p\xi\right),$$
$$h_{P}(T,\xi)={\limsup_j} \ h_j(T,\xi),$$
$$h_{P}(T)=\sup_\xi h_{P}(T,\xi),$$
where $\xi$ denotes a finite measurable partition of the probability space, and $H(\xi)$ is the entropy of the partition $\xi$:
$$ H(\{C_1,C_2,\dots, C_n\})=-\sum_{i=1}^n \mu( C_i)\ln \mu( C_i).$$

We will consider a special case of Kirillov-Kushnirenko entropy: the case of the action of the group $Z$, and the sets $P_j$ are chosen in the form of increasing progressions
 $P_j=\{j,2j,\dots, L(j)j\},$ $ L(j)\to\infty.$

An automorphism $G$ has \bf completely positive P-entropy \rm if
$h_{P}(G,\xi)>0$ for any non-trivial partition $\xi$.

\vspace{3mm}
\bf Theorem 1. \it If an automorphism $S$ has zero entropy, then there is a Gaussian automorphism $G$ and the sequence
$P_j=\{j,2j,\dots, L(j)j\},$ $ L(j)\to\infty,$ such that
$h_P(S)=0$ and $G$ has completely positive $P$-entropy.
Such $S$ and $G$ are disjoint: they do not admit nontrivial Markov wreath products. \rm

\vspace{3mm}
\bf Lemma (\cite{21}). \it For a deterministic $S$ there is a sequence of the form $P_j=\{j,2j,\dots, L(j)j\},$ $ L(j)\to\infty,$ such that $h_P(S) =0$. \rm

\vspace{3mm}
The disjointness of $S$ and $G$ follows from Theorem 4.1 of \cite{RT}.
So we need to find a suitable $G$. Note that in \cite{RT}
a suitable Poisson superstructure $P(T)$ was found over a specially selected transformation $T:\R\to\R$ of rank 1. If we consider $T$ as an orthogonal operator on the space $L_2(\R)$, which is endowed with a Gaussian measure, then the operator $T$ induces an action $G(T)$ that preserves this measure. 

Let an automorphism $S$ be of  zero entropy. By virtue of  Lemma, we find a suitable sequence $P'=P_{j'}$ such that $h_{P'}(S)=0$.

 Following \cite{RT} and \cite{24} for the sequence $P_{j'}$ we can 
  find (see below)  an orthogonal operator $T$ and a subsequence $P=\{P_{j}\}$ such that the following conditions are satisfied:

(i)$T$ has a singular spectrum,

(ii) the spectrum $T\otimes T$ is absolutely continuous,

(iii) for $T$ the condition $(\perp^P)$ is satisfied.

Definition. Let $H_j$ be a sequence of expanding finite-dimensional spaces in real $L_2(X)$ whose union is dense in $L_2(X)$. Let the orthogonal operator $T$ satisfy

($\perp^P$): \it spaces $T^pH_j$, $p\in P_j$, are pairwise orthogonal. \rm

Such $T$ exists among Sidon transformations. Properties
(i),(ii) have the constructions from \cite{24}, they are compatible
with property (iii), which appeared in implicit form in the work of \cite{RT}. Let us recall what was made there. Some transformation $T$ acted on $X=\cup_j X_j$, where $X_j$ is an increasing sequence of sets whose measures tend to infinity.
For each $j$, the sets $T^pX_j$, $p\in P_j$, do not intersect.
Therefore, the spaces $T^pL_2(X_j)$, $p\in P_j$, where $L_2(X_j):=\{f\chi_{X_j}:\, f\in L_2(X)\}$ are pairwise orthogonal.
We choose the required finite-dimensional $H_j$ as subspaces of the spaces $L_2(X_j)$.

For a set $A\subset H_j$ measurable with respect to the Lebesgue measure on $H_j$, we define the cylindrical set $C_A=\{f:\, \pi_j f \in A\}$, where $\pi_j $ is the orthoprojection on $ H_j$.
Due to the indicated orthogonality, we note that all sets
$T^pC_A$, $p\in P_j$, are independent with respect to the Gaussian measure.
The action of $T$ on cylinders $C_A$ is now denoted by $G(T)$.
Due to the noted independence and asymptotical density of $H_j$, we obtain that $G(T)$ has  completely positive P-entropy.

Thus, we can present an infinite set (in fact, of continual cardinality) of pairwise disjoint Gaussian automorphisms with  spectrally similar even factors (with Lebesgue spectrum). 

\section{Гауссовские автоморфизмы с вполне положительной P-энтропией}
\bf $P$-энтропия (случай энтропии Кириллова-Кушниренко). \rm 
Пусть $P=\{P_j\}$ — последовательность конечных подмножеств в счетной бесконечной группе. Для сохраняющего меру действия $T$ этой группы на вероятностном пространстве определим
$$h_j(T,\xi)=\frac 1 {|P_j|} H\left(\bigvee_{p\in P_j}T_p\xi\right),$$
$$h_{P}(T,\xi)={\limsup_j} \ h_j(T,\xi),$$
$$h_{P}(T)=\sup_\xi h_{P}(T,\xi),$$
где $\xi$ обозначает конечное измеримое разбиение вероятностного пространства, а $H(\xi)$ — энтропия разбиения $\xi$:
$$ H(\{C_1,C_2,\dots, C_n\})=-\sum_{i=1}^n \mu( C_i)\ln \mu( C_i).$$

Мы   рассмотрим  случай действия группы $Z$, причем  множества  $P_j$ выбираются в виде  возрастающих прогрессий
    $$P_j=\{j,2j,\dots, L(j)j\}, \ \ \ L(j)\to\infty.$$

Автоморфизм $G$ имет вполне положительную энтропию, если 
$h_{P}(G,\xi)>0$ для всякого нетривиального разбиения $\xi$.

\vspace{3mm}
\bf Теорема 1. \it Если автоморфизм $S$ имеет нулевую энтропию, то найдется гауссовский автоморфизм $G$ и последовательность 
$P_j=\{j,2j,\dots, L(j)j\},$ $ L(j)\to\infty,$ такие, что 
$h_P(S)=0$ и $G$ обладает вполне положительной $P$-энтропией. 
Такие $S$ и $G$ дизъюнктны: они не допускают нетривиальные марковские сплетения.   \rm

\vspace{3mm}
\bf  Лемма (\cite{21}). \it  Для детерминированного $S$ найдется последовательность вида  $P_j=\{j,2j,\dots, L(j)j\},$ $ L(j)\to\infty,$  такая, что $h_P(S)=0$.   \rm

\vspace{3mm}
Дизъюнктность $S$ и $G$ вытекает из  теоремы 4.1 работы \cite{RT}.
Таким образом нам нужно найти подходящий $G$.  Отметим, что в \cite{RT}
была найдена подходящая пуассоновская надстройка $P(T)$ над специально подобранным преобоазованием $T:\R\to\R$ ранга 1. Если же рассматривать $T$  как ортогональный оператор на пространстве $L_2(\R)$, которое наделено гауссовской мерой, то  оператор $T$ индуцирует действие $G(T)$, которое сохраняет эту меру.  Энтропийные свойства $P(T)$ и $G(T)$
 в интересующей нас ситуации оказываются похожими.

Итак, дан автоморфизм $S$  с нулевой энтропией.   В силу леммы находим подходящую последовательность $P'=P_{j'}$, чтобы $h_{P'}(S)=0$.

 Cледуя \cite{RT} для    последовательности $P_{j'}$  найдем ортогональный  оператор $T$ и  подпоследовательность $P=\{P_{j}\}$ такие, что выполнены условия:  

(i)$T$ обладает сингулярным  спектром,

(ii)  спектр $T\otimes T$ абсолютно непрерывен,

(iii)  для  $T$ выполнено условие $(\perp^P)$. 

Определение.  Пусть  $H_j$ -- последовательнось расширяющихся конечномерных пространств в $L_2(X)$, объединение которых плотно в $L_2(X)$.  Пусть для ортргонального оператора $T$ выполнено 

($\perp^P$): \it пространства $T^pH_j$, $p\in P_j$, попарно ортогональны. \rm

Такой $T$  имется среди сидоновских преобразований. Свойствами 
(i),(ii) обладают конструкции из работы \cite{24}, они совместимы
со свойством (iii), которое в неявной форме  фигурировало в работе \cite{RT}.
Напомним, о чем там шла речь.  Некоторое преобразование $T$  действовало на $X=\cup_j X_j$, где $X_j$ -- возрастающая последовательность множеств,  меры которых стремятся к бесконечности.
Для каждого $j$ множества $T^pX_j$, $p\in P_j$, не пересекаются.
Поэтому пространства $T^pL_2(X_j)$, $p\in P_j$, где  $L_2(X_j):=\{f\chi_{X_j}:\, f\in L_2(X)\}$ попарно ортогональны.
Нужные конечномерные  $H_j$ выбираем как подпространства пртранств $L_2(X_j)$.

Для множества $A\subset H_j$,  измеримого относительно меры Лебега на $H_j$, определим  цилиндрическое множество $C_A=\{f:\, \pi_j f \in A\}$, где $\pi_j $-- ортопроекция на  $H_j$.
В силу указанной ортогональности замечаем, что все множества 
$T^pC_A$, $p\in P_j$, независимы отностительно гауссовской меры.
Действие $T$ на цилиндрах $C_A$ теперь обозначается через $G(T)$.
В силу отмеченной независимости получается, что $G(T)$ обладает 
вполне положительной P-энтропией.

Таким образом, мы можем предъявить бесконечно много (на самом деле континуум) попарно неизоморфных, но спектрально одинаковых  четных факторов гауссовских автоморфизмов. Для пуассоновских надстроек ситуация иная.  Метрические свойства $G(T)$ и $P(T)$ могут сильно различаться: например, гауссовы автоморфизмы бесконечно делимы, следовательно, $G(T)$ имеет обширный  централизатор, а центраизатор $P(T)$ для наших $T$ будет  тривиальным (см. \cite{06}, \cite{15},\cite{PR}).

\large

\end{document}